\documentclass[makeidx]{amsart}

\usepackage{a4wide}
\usepackage{amssymb,amscd}
\usepackage[utf8]{inputenc}
\usepackage{pb-diagram}
\usepackage{graphicx}
\usepackage{color}
\usepackage[hyperref]{hyperref}

\renewcommand{\phi}{\varphi}
\newcommand{\C}{{\mathbb{C}}}
\newcommand{\N}{{\mathbb{N}}}
\newcommand{\R}{{\mathbb{R}}}

\newcommand{\Z}{{\mathbb{Z}}}
\newcommand{\1}{{\mathbf{1}}}

\renewcommand{\epsilon}{\varepsilon}
\renewcommand{\theta}{\vartheta}
\renewcommand{\S}{{\mathbb{S}}}

\newcommand{\CP}[1]{{\C\mathbb{P}(#1)}}
\newcommand{\invOrientCP}[1]{{\overline{\C\mathbb{P}}(#1)}}
\newcommand{\Reeb}{{X_\mathrm{Reeb}}}
\newcommand{\lie}[1]{{\mathcal{L}_{#1}}}

\newcommand{\norm}[1]{{\lVert #1\rVert}}
\newcommand{\abs}[1]{{\left\lvert #1\right\rvert}}

\newcommand{\z}{{\mathbf{z}}}

\DeclareMathOperator{\sing}{Sing}

\DeclareMathOperator{\stab}{Stab}

\theoremstyle{plain}
\newtheorem{theorem}{Theorem}
\newtheorem{propo}[theorem]{Proposition}
\newtheorem{coro}[theorem]{Corollary}

\theoremstyle{remark}
\newtheorem{remark}{Remark}
\newtheorem{example}{Example}

\theoremstyle{definition}
\newtheorem*{defi}{Definition}

\numberwithin{equation}{section}

\begin{document}

\title{Resolution of symplectic cyclic orbifold singularities}

\author[K.\ Niederkrüger]{Klaus Niederkrüger}
\author[F.\ Pasquotto]{Federica Pasquotto}

\email[K.\ Niederkrüger]{kniederk@umpa.ens-lyon.fr}
\email[F.\ Pasquotto]{pasquott@few.vu.nl}

\address[K.\ Niederkrüger]{École Normale Supérieure de Lyon\\
  Unité de Mathématiques Pures et Appliquées\\
  UMR CNRS 5669\\
  46, allée d'Italie\\
  69364 LYON Cedex 07\\
  France}

\address[F.\ Pasquotto]{Department of Mathematics, Faculty of Sciences\\
  Vrije Universiteit\\De Boelelaan 1081a\\1081 HV Amsterdam\\
  The Netherlands}

\begin{abstract}
  In this paper we present a method to obtain resolutions of
  symplectic orbifolds arising from symplectic reduction of a
  Hamiltonian $\S^1$--manifold at a regular value.
  
  As an application, we show that all isolated cyclic singularities of
  a symplectic orbifold admit a resolution and that pre-quantisations
  of symplectic orbifolds are symplectically fillable by a smooth
  manifold.
\end{abstract}


\maketitle

\section{Introduction}

Symplectic quotients are an important source of new symplectic
manifolds: they appear as \emph{symplectic reductions} in the context
of Hamiltonian actions and the associated moment maps.  More
generally, reduced spaces corresponding to regular values of the
moment map turn out to be \emph{symplectic orbifolds}, but one can
still look for a closed smooth symplectic manifold which is isomorphic
to the orbifold outside a neighbourhood of the singular set: we call
this object a \emph{symplectic orbifold resolution}.

Even in the case of a reduced space corresponding to a singular value
of the moment map, where singularities of a more complicated type can
occur, Kirwan's ``partial desingularisation'' method
\cite{KirwanDesingularisation} can be applied to obtain a resolution
which has only orbifold singularities.

The method described in this paper relies on the construction of an
auxiliary circle action in a neighbourhood of the singularities with
largest structure group: this action is subsequently used to perform a
symplectic cut (as described in \cite{Lerman_SymplecticCuts}), which
also amounts to a weighted blow up along the singular stratum
\cite{Godinho_blowup}.  The orbifold obtained in this way has
singularities of lower order, and we can repeat the step inductively
until we produce a smooth manifold. Making use of this procedure, we
are able to find orbifold resolutions for reduced spaces obtained from
symplectic reduction at the regular level sets of a Hamiltonian
function generating an $\S^1$--action.  In particular, we find
symplectic resolutions for all isolated cyclic orbifold singularities.
A technique involving symplectic resolutions of orbifolds via blow-ups
can also be found in \cite{Cavalcanti_Fernandez_Munoz_NonFormal}: due
to a very restrictive definition of orbifolds, though, that result
only provides resolutions in the special case of isolated
singularities.

As a direct application of our result, we are able to show that a
Seifert manifold with a contact structure that is $\S^1$--invariant
and transverse to the fibres is symplectically fillable by a smooth
manifold.

\subsubsection*{Acknowledgements}

K.\ Niederkrüger is working at the \emph{ENS de Lyon}, where he is
being funded by the Agence Nationale de la Recherche (ANR) project
\emph{Symplexe}.  F.\ Pasquotto is working at the \emph{Vrije
  Universiteit Amsterdam} and is supported by VENI grant 639.031.620
of the Nederlandse Organisatie voor Wetenschappelijk Onderzoek (NWO).

We thank the \emph{Mathematisches Forschungsinstitut Oberwolfach} for
inviting us to a \emph{Research in Pairs} stay at the beginning of
2007, during which we developed the main parts of this article.

\subsection{Isolated singularities in dimension four}

We start our paper by discussing singularities of the simplest form:
we hope that this will provide the reader with some motivation and
will serve as the right introduction to the difficulties arising when
considering more general examples.

In the four-dimensional case, we can give a very explicit description
of the resolution of an isolated orbifold singularity.  In order to do
so, we use weighted blow-ups of isolated symplectic orbifold
singularities, as defined in \cite{Godinho_blowup}.

Let $x\in \left(M^{(4)},\omega\right)$ be an isolated orbifold
singularity with structure group $\Z_p=\{ 1, \xi, \xi^2,\dotsc ,
\xi^{p-1}\}$.  The orbifold chart around $x$ is $\Z_p$--equivariantly
symplectomorphic to a neighbourhood of $0$ in $\C^2$ with $\Z_p$
acting by
\begin{equation*}
  (z_1, z_2)\mapsto (\xi^m z_1, \xi z_2),\quad 0<m<p, \quad \gcd(m, p)=1\;.
\end{equation*}
If $m$ and $p$ were not coprime, the orbifold singularity would not be
isolated.  We can define an $\S^1$--action on $\C^2$ by setting
$\lambda\cdot (z_1, z_2)=(\lambda^m z_1, \lambda z_2)$ for
$\lambda\in\S^1$.  The actions of $\Z_p$ and $\S^1$ commute, but the
induced circle action on $\C^2/\Z_p$ is not effective.  Instead we
have to go to the $\S^1/\Z_p$--action obtained from the following
exact sequence
\begin{equation*}
  0\longrightarrow \Z_p\longrightarrow\S^1\longrightarrow \hat\S^1
  \longrightarrow 0 \;,
\end{equation*}
with the homomorphism of the circle given by $\lambda\mapsto
\lambda^p$.  This defines now a symplectic $\hat\S^1$--action on
$\C^2/\Z_p$ by
\begin{equation*}
  \mu \cdot[z_1, z_2]=[\lambda\cdot (z_1, z_2)]=[\lambda^m z_1,
  \lambda z_2]
\end{equation*}
for $\mu\in\hat\S^1$ and a $\lambda\in \S^1$ such that
$\lambda^p=\mu$.

The weighted blow-up of $\C^2/\Z_p$ at the origin can be represented
as a symplectic cut with respect to this $\hat\S^1$--action: Take the
product manifold $\C^2/\Z_p \times \C$ with the symplectic form given
by $(\omega, -i\,dw\wedge d\bar w)$ and the effective
$\hat\S^1$--action
\begin{equation*}
  \mu \cdot\bigl([z_1, z_2], w\bigr) = \bigl(\mu\cdot [z_1, z_2],
  \mu^{-1} w\bigr) = \bigl([\lambda^m z_1, \lambda z_2], \lambda^{-p}
  w\bigr), \quad \lambda^p=\mu \;.
\end{equation*}
The blow-up is the symplectic reduction of this space.  The
Hamiltonian function that corresponds to this action is
\begin{equation*}
  H\bigl([z_1,z_2],w\bigr) = m\abs{z_1}^2 + \abs{z_2}^2 - p\abs{w}^2\;,
\end{equation*}
and the $\epsilon$--level set $H^{-1}(\epsilon)$ is diffeomorphic to
the manifold $\S^3/\Z_p\times \C$.  The last step consists in taking
the quotient $H^{-1}(\epsilon)/\hat\S^1$ to obtain a symplectic
orbifold that can be glued to $\C^2/\Z_p$ after removing a
neighbourhood of the origin.

There is only one singular point in $H^{-1}(\epsilon)/\hat\S^1$,
namely $\bigl([1, 0], 0\bigr)\in \S^3/\Z_p\times \C$, with stabiliser
$\Z_m$.  By the slice theorem, a neighbourhood of this point admits an
orbifold chart equivalent to $\C^2$ with structure group $\Z_m$ acting
by $\eta\cdot (w_1, w_2)=(\eta^{-p} w_1, \eta w_2)$.

Now choose $a_1\in \Z$ such that $0<a_1 m -p<m$ and set $m_1 = a_1 m -
p$ and $p_1 = m$: then the new singularity can also be modelled by
$\Z_{p_1}$ acting by $\eta\cdot (z_1, z_2)= (\eta^{m_1} z_1, \eta
z_2)$ with $\gcd(m_1,p_1) = 1$, because if $b$ divides both $p_1$ and
$m_1$ then it also divides $m$ and $p$.  We are thus in the initial
type of situation, but we have managed to reduce the order of the
singularity.

If we blow up once more, we replace this by a new singularity, this
time with structure group $\Z_{m_1}$ acting by $\zeta\cdot (z_1,
z_2)=\ \bigl(\zeta^{a_2m_1 -p_1} z_1, \zeta z_2\bigr)$.  If we iterate
this blow-up process, at each step we replace the previous singularity
by a new one with structure group $\Z_{p_i}$ acting by $(\xi^{m_i}
z_1, z_2)$, where the pair $(p_i, m_i)$ is recursively given by
\begin{equation*}
  \begin{pmatrix}
    p_i \\ m_i
  \end{pmatrix}
  =
  \begin{pmatrix}
    0 & 1 \\ -1 & a_i
  \end{pmatrix}
  \begin{pmatrix}
    p_{i-1} \\ m_{i-1}
  \end{pmatrix}
  \;,
\end{equation*}
with each $a_i$ corresponding to the ``roundup'' of
$\frac{p_{i-1}}{m_{i-1}}$, that is, the least integer $\geq
\frac{p_{i-1}}{m_{i-1}}$. The sequence $[a_1, a_2, \dotsc]$
corresponds to the continued fraction of $\frac{p}{m}$, so in
particular it is a finite sequence (a description of resolutions in
terms of continued fractions is contained, for example, in Miles
Reid's lecture notes \cite{ReidLectureNotes}).  After sufficiently
many blow-ups, in other words, we get a pair of the form $(p_N, 1)$
and thus an orbifold chart which is in fact smooth. This is then the
resolution of our initial singularity.
\begin{figure}[htbp]
  \begin{center}
    \input{cyclic_dim4_resolution.pstex_t}
    \caption{The resolution obtained via a sequence of blow-up can
      also be thought of as the connected sum $ M \#
      \invOrientCP{1,m,p} \# \invOrientCP{1,m_1,p_1} \# \dotsm \#
      \invOrientCP{1,1,p_N}$.}  \label{fig: resolution_4dim}
  \end{center}
\end{figure}

Notice that we can think of each weighted blow-up as taking the
connected sum with a suitable orbifold. In dimension $4$, if we define
the \emph{weighted projective space} to be
\begin{equation*}
  \CP{a_0, a_1, a_2} = \Bigl\{[z_0:z_1:z_2]\sim
  [\lambda^{a_0} z_0:\lambda^{a_1} z_1: \lambda^{a_2} z_2]\,\Bigm|\,
  (z_0, z_1, z_2)\in \C^3-\{0\}\text{ and }\lambda\in \C^*\Bigr\},
\end{equation*}
the weighted blow-up of a singular point with structure group $\Z_p$
acting by $(z_1, z_2)\mapsto (\xi^a z_1, \xi^b z_2)$ can be described
as taking the connected sum, around this point, with the orbifold
$\CP{a,b,p}$ with reversed orientation.  We can use this description
to represent the resolution of a four-dimensional cyclic singularity
as in Fig.~\ref{fig: resolution_4dim}.

In higher dimension ($\geq 6$), the method just described is not
suitable: even if we start with an isolated singularity, after the
first blow-up the singular set is not necessarily discrete any more.

\begin{example}
Consider an isolated orbifold singularity modelled on a neighbourhood 
of $0$ in $\C^3$ with $\Z_p$ acting by
\begin{equation*}
  (z_1, z_2, z_3)\mapsto (\xi^{m_1} z_1, \xi^{m_2} z_2, \xi z_3),
\quad 0<m_1<m_2<p, \quad \gcd(m_j, p)=1, j=1,2\;.
\end{equation*}
We can blow up this singularity with the method used in the
four-dimensional case, that is, performing a symplectic cut with
respect to a suitable circle action.  If $\gcd(m_1, m_2)=d\neq 1$,
though, after blowing up the singular set will be a two-dimensional
suborbifold with stabiliser $\Z_d$ at generic points.
\end{example}

One could still hope to obtain a resolution of general cyclic orbifold
singularities by using weighted blow-ups along suborbifolds
(cf.~\cite{Meinrenken_Sjamaar}), but for this we would need to find a
suitable circle action on the fibres of the normal orbibundle to
singular strata.  While this is not difficult in individual bundle
charts, we were not able to define a global action.

In the rest of this paper we thus adopt a different desingularisation
method, which has the advantage of applying to a larger class of
orbifold singularities than just cyclic isolated ones.

\section{Symplectic orbifolds as quotients of Hamiltonian
  $\S^1$--manifolds}

Let $M$ be a manifold with an action of the unit circle $\S^1$.  The
infinitesimal action is defined by the vector field
\begin{equation*}
  X_M(p) := \left.\frac{d}{dt}\right|_{t=0} \exp(tX)*p
\end{equation*}
for every $p\in M$.

\begin{defi}
  A \textbf{Hamiltonian $\S^1$--manifold} $(M,\omega)$ is a symplectic
  manifold with an $\S^1$--action for which there exists a Hamiltonian
  function $H:\, M\to\R$ that generates the action or, equivalently,
  that satisfies the identity
  \begin{equation*}
    i_{X_M}\omega = dH \;.
  \end{equation*}
  If this is the case, the function $H$ is also called the
  \textbf{moment map} of the action.
\end{defi}

The following statement is well known (\cite{WeinsteinVmanifolds}),
but for completeness we will briefly sketch the proof.

\begin{propo}
  Let $(M,\omega)$ be a Hamiltonian $\S^1$--manifold with Hamiltonian
  function $H$.  For any regular value $E$ of $H$, the
  symplectic quotient
  \begin{equation*}
    M_E := H^{-1}(E) / \S^1
  \end{equation*}
  is a symplectic orbifold.
\end{propo}
\begin{proof}
  Since $E$ is a regular value, it follows that $H^{-1}(E)$ is a
  smooth submanifold on which the $\S^1$--action is semi-free, because
  $i_{X_M} \omega = dH \ne 0$.  The stabilizer of a point $p\in
  H^{-1}(E)$ is thus isomorphic to some finite cyclic subgroup.  By
  the slice theorem, we can find a homeomorphism from a neighbourhood
  of the equivalence class of $p$ in the orbit space to the quotient
  $S_p/\stab(p)$, where $S_p$ denotes a slice at $p$ (the normal space
  to the $\S^1$--orbit). In this way we obtain an orbifold atlas for
  the space $M_E$.  Moreover, the symplectic form $\omega$ induces a
  symplectic structure on the orbifold $M_E$, because $\lie{X_M}\omega
  \equiv 0$, and $\left.  \omega(X_M,\cdot)\right|_{TH^{-1}(E)} =
  \left.  dH\right|_{TH^{-1}(E)} \equiv 0$.
\end{proof}

\begin{example}
  The \textbf{weighted projective space} $\CP{a_0, a_1,\dotsc , a_n}$
  is the $2n$--dimensional orbifold obtained as the symplectic
  quotient of $\C^{n+1}$ by the $\S^1$--action
  \begin{equation*}
    \lambda*(z_0, z_1,\dotsc,z_n)=(\lambda^{a_0} z_0,\dotsc,
    \lambda^{a_n} z_n) \;.
  \end{equation*}
  It is in fact a symplectic orbifold with the structure induced by
  the canonical symplectic form on $\C^{n+1}$, namely $\omega_0=
  \frac{i}{2}\sum dz_j\wedge d\bar{z}_j$. The Hamiltonian function
  generating the $\S^1$--action is
  $H(z_0,\dotsc, z_n) = a_0\abs{z_0}^2 + \dotsm + a_n\abs{z_n}^2$.
\end{example}

\subsection{Symplectic resolutions} 

\begin{defi}
  Let $(M, \omega)$ be a symplectic orbifold with singular set $X$.  A
  \textbf{symplectic $\epsilon$--resolution} of $M$ consists of a
  smooth symplectic manifold $(\widetilde{M}, \widetilde{\omega})$ and
  a continuous map $p:\, (\widetilde{M}, \widetilde{\omega})\to (M,
  \omega)$ for which there is an orbifold tubular neighbourhood $U$ of
  $X$ of size $\epsilon$, such that
  \begin{equation*}
    \left.p\right|_{\widetilde{M}-p^{-1}(U)}:\,
    \widetilde{M}-p^{-1}(U)\to M-U
  \end{equation*}
  is a symplectic diffeomorphism.
\end{defi}

In this paper we prove the following result.

\begin{theorem}\label{main}
  Let $M$ be a symplectic orbifold arising from symplectic reduction
 of a Hamiltonian $\S^1$--manifold at a regular value of the moment map.
  Then there exists a symplectic $\epsilon$--resolution of~$M$ for any
  arbitrarily small~$\epsilon>0$.
\end{theorem}

\subsection{The stratification of the singular set}

Let $W$ be a symplectic $\S^1$--manifold with Hamiltonian $H$ and
denote by $\sing W\subset W$ the singular set of the circle action,
that is,
\begin{equation*}
  \sing W = \bigl\{x\in W\,\bigm|\,\stab(x)\neq \{1\}\bigr\} \;.
\end{equation*}
For a given isotropy group $\Z_k$, let $W_k$ denote the union of
orbits whose isotropy group is $\Z_k$, namely $W_k=\{x\in
W\,|\,\stab(x)\cong \Z_k\}$.  Then $\sing W$ is stratified by singular
strata, i.e., connected components of $W_k$ for all $k\neq 1$.  Assume
that $0$ is a regular value of $H$ and denote by $P$ the level set
$H^{-1}(0)$ and by $\sing P$ the singular set of the action on this
level set.  If $\pi:\,P\to M = P/\S^1$ denotes the orbit map,
$X=\pi(\sing P)$ is the set of orbifold singularities of $M$ and the
stratification of $\sing P$ descends to a stratification of~$X$.  Our
desingularisation method works by induction on the order of the
stabilisers of the strata $W_k$.  Namely, we start with the stratum
with largest stabiliser (\emph{minimal stratum}), remove a small
neighbourhood of it and glue in a smooth manifold, in such a way that
both the symplectic form and the Hamiltonian $\S^1$--action extend to
the manifold resulting from this surgery.  The new singular set will
also carry a stratification, but the order of the stabiliser of the
minimal stratum will be strictly less than $k$.  If we successively
repeat this procedure sufficiently often, we will eventually reduce
this maximal stabiliser to the trivial group. In particular,
symplectic reduction at $0$ yields then a smooth symplectic manifold
which gives a resolution of the orbifold~$M$.

After working out this desingularisation method for symplectic
orbifolds, we were told by one of the authors of
\cite{GinzburgGuilleminKarshon} that they, and others before them, had
already used the same construction in the smooth category.

\section{Construction of the resolution}

The strategy of our construction is to find an auxilliary circle
action on the Hamiltonian $\S^1$--manifold around its minimal stratum.
This action will have certain properties (see
Proposition~\ref{properties of phi-action}) that will allow us to
perform a symplectic cut \cite{Lerman_SymplecticCuts}. The stabilizers
of the Hamiltonian $\S^1$--manifold obtained this way are all lower
than the stabilizer of the initial minimal stratum.  Then we
successively repeat the construction until we obtain a manifold with a
free circle action.

\subsection{An auxiliary circle action around the minimal stratum}
\label{sec: construction of auxiliary circle action}

Suppose $(W,\omega)$ is a Hamiltonian $\S^1$--manifold with
Hamiltonian function $H:\,W\to \R$.  Assume further that $0$ is a
regular level set of $H$, and choose a metric and an almost complex
structure $J$ which are $\S^1$--invariant and compatible with $\omega$
(see for example \cite[Section~5.5]{McDuffSalamonIntro}).  Since
symplectic reduction is a local process, it will be sufficient to 
study a neighbourhood of $P:=H^{-1}(0)$ consisting of regular level
sets $P_t:=H^{-1}(t)$, $t\in (-\epsilon, \epsilon)$.  We will start by
finding a suitable model for this neighbourhood.

\begin{propo}
  There is a neighbourhood of $P$ in $(W,\omega)$ that is
  $\S^1$--diffeomorphic to
  \begin{equation*}
    P\times (-\epsilon,\epsilon) \;,
  \end{equation*}
  with trivial circle action on the second factor $\lambda *(p,t) =
  (\lambda *p,t)$. Under this diffeomorphism, the Hamiltonian function
  pulls back to $H(p,t) = t$.
\end{propo}
\begin{proof}
Consider the vector field
  \begin{equation*}
    Y := \frac{1}{\norm{\nabla H}^2}\, \nabla H
    = \frac{1}{\norm{X_{\S^1}}^2}\, JX_{\S^1}\;.
  \end{equation*}
  It is transverse to the level sets, and allows us to define for
  small times $t$ with $\abs{t}<\epsilon$ the diffeomorphism
  \begin{equation*}
    \Psi:\,P\times(-\epsilon,\epsilon) \to W, \, (p,t) \mapsto \Phi^Y_t(p) 
  \end{equation*}
  onto a neighbourhood of $P$ in $W$.  The pullback of $H$ gives $H\circ
  \Psi(p,t) = t$, and the circle action is $t$--invariant (by which we
  mean $\lambda *\Psi(p,t) = \Psi(\lambda*p,t)$ for $\lambda
  \in\S^1$), because $Y$ commutes with the generator of the
  $\S^1$--action.
\end{proof}

>From now on, we will always assume the neighbourhood of the level set
to be of the form $P\times (-\epsilon,\epsilon)$ with an
$\S^1$--action on the first factor, and with Hamiltonian function
$H(p,t) = t$.

In the next step, we will construct a suitable model for a
neighbourhood of the minimal stratum, in which we can describe the
auxiliary circle action needed for the resolution.

The $\S^1$--action on $W$ is semi-free.  Consider the stratification
$\{W_k\}$ of the singular set $\sing W$ given by the isotropy groups.
Each stratum $W_k$ is an $\S^1$--invariant symplectic submanifold of
$W$ of the form $P_k\times (-\epsilon,\epsilon)$, where $P_k :=
W_k\cap P$.

If $k$ is maximal, that is, there are no points in $W$ of order larger
than $k$, then $W_k\cap P_t$ is a closed submanifold. By restricting
to one component, we may further assume $W_k$ to be connected.  Choose
any $\S^1$--invariant metric $g_P$ on $P$ and extend it to the product
metric $g_W := g_P\oplus dt^2$ on $W = P\times (-\epsilon,\epsilon)$.
With our choice of metrics on $P$ and $W$, the corresponding
exponential maps are related as follows:
\begin{equation*}
  \exp_{(p,t)}^W (v+a\,\partial_t) = \bigl(\exp_p^P (v), t+a \bigr)
\end{equation*}
with respect to the splitting $T_{(p,t)} W = T_p P \oplus \langle
\partial_t\rangle$.  Moreover, if we denote the normal bundle of $P_k$
in $P$ by $N_k$, then the normal bundle of $W_k$ in $W$ is
isomorphic to the product $N_k\times (-\epsilon, \epsilon)$.
Therefore a tubular neighbourhood of $W_k$ in $W$ can also be
identified via the exponential map with the product $N_k(\delta)
\times (-\epsilon, \epsilon)$, where $N_k(\delta)$ denotes the
$\delta$--subdiskbundle of $N_k$.  This can be summarized in the
following proposition:

\begin{propo}
  There exists a neighbourhood of the singular stratum $W_k$ in $W$
  that has the form $N_k(\delta)\times (-\epsilon,\epsilon)$.  The
  Hamiltonian function on this neighbourhood is just given by $H(v,t)
  = t$ and the circle action is the linearized $\S^1$--action
  on~$N_k$.
\end{propo}

In what follows we will therefore implicitly assume this
identification and denote for simplicity by $\omega$ and $J$ the
pullback of the corresponding structures on $W$.

Later it will be necessary to introduce a second circle action.  To
avoid confusions, from now on we will call the linearisation of the
given action the $\beta$--action, and we will write it as $\lambda
*_\beta v$ for $\lambda\in\S^1$ and~$v\in \nu_k(\delta)$.

Let $x\in W_k$ be a singular point in the minimal stratum.  The
stabiliser $\stab(x)\cong\Z_k$ acts by isometric (with respect to both
$g_W$ and the metric $g_J := \omega(\cdot, J\cdot)$), $J$--linear
transformations on $T_x W$, hence the $\Z_k$--action is equivalent to
\begin{equation*}
  \lambda*_\beta \z = (\lambda^{\widetilde{a}_1}z_1,\dotsc,
  \lambda^{\widetilde{a}_n}z_n), 
  \quad \lambda\in\Z_k, \; \z\in T_xW,
  \;\widetilde{a}_1,\dotsc,\widetilde{a}_n\in\Z\;.
\end{equation*}
Without loss of generality we may assume that $0 = \widetilde{a}_1 =
\dotsb = \widetilde{a}_m < \widetilde{a}_{m+1}\leq \dotsb \le
\widetilde{a}_n < k$ for some $m$.  The first $m$ directions span the
space $T_xW_k$, and the other directions are orthogonal to $T_x W_k$:
it is easy to show that this holds not only for $g_J$, but also with
respect to any other $\S^1$--invariant metric. Hence, if we denote by
$\nu_k$ the normal bundle of $W_k$ in $W$, the subspace $\{(0,\dotsc ,
0,z_{m+1}, \dotsc, z_n)\}$ coincides with the fibre $\nu_k(x)$.
For all $x\in W_k$, $j_x := \left.J\right|_{\nu_k(x)}$ defines a
complex vector bundle structure on $\nu_k$, making the $\Z_k$--action
$j$--complex linear in each fibre.

Denote by $a_1<\dotsb <a_l$ the \emph{distinct} exponents occurring in
the normal form for the action: $\nu_k$ splits thus into a direct sum
of subbundles
\begin{equation*}
  \nu_k = E_1\oplus\dotsb \oplus E_l \;,
\end{equation*}
where $E_i(x)$ denotes the eigenspace corresponding to the eigenvalue
$\lambda^{a_i}$ in the fibre at the point $x$.  This splitting is well
defined for each component of~$W_k$.  Therefore we can extend the
$\Z_k$--action to a second circle action by setting for any $\lambda
\in\S^1$
\begin{equation*}
  \lambda *_{\mathrm{aux}} v  := \lambda^{a_1}v_1+ \dotsb +\lambda^{a_l}v_l \;,
\end{equation*}
where $v=v_1+\dotsb +v_l$ is a splitting with respect to the
eigenspaces defined above, and $\lambda^{a_i}v_i := \bigl(\cos
(a_i\lambda) + \sin (a_i\lambda)\,j\bigr) v_i$.  This auxiliary action
is fibrewise and $j$--linear, and commutes with the original
$\beta$--action.  Unfortunately it does not need to respect the
symplectic form~$\omega$.  Recall that $\omega$ denotes here the
symplectic form on the total space of $\nu_k$ obtained by pulling back
the symplectic form on $W$.  By averaging $\omega$ over the auxiliary
action, we obtain a closed $2$--form $\omega^\prime$ on $\nu_k$ that
is invariant with respect to both the $\beta$-- and the auxiliary
action.  There is a small neighbourhood of the zero section of $\nu_k$
where we also have $(\omega^\prime)^n\neq 0$ and hence $\omega^\prime$
is a symplectic form.  To prove this it suffices to show that $\omega$
is ${\mathrm{aux}}$--invariant on the zero section.  In fact, at the
zero section $W_k$ of $\nu_k$, there is a well defined splitting of
the tangent bundle $\left.T\nu_k\right|_{W_k} = T W_k \oplus \nu_k$
and we can write
\begin{equation*}
  \left.T\nu_k\right|_{W_k} = T W_k \oplus E_1\oplus\dotsb \oplus
  E_l\;.
\end{equation*}
The $E_i$'s are $J$--linear subspaces and $g_J$--orthogonal to each
other, so that they are also symplectically orthogonal.  The
linearised auxiliary action on a vector $w+v_1+\dotsm + v_l\in T\nu_k$
with $w\in T_x W_k$ and $v_i\in E_i(x)$ is given by
\begin{equation*}
  \lambda *_{\mathrm{aux}} (w+v_1+\dotsm + v_l) = w + \lambda^{a_1}v_1+ \dotsb
  +\lambda^{a_l}v_l
\end{equation*}
and using the orthogonality relations, it is easy to check that
\begin{equation*}
  \omega (\lambda *_{\mathrm{aux}} v, \lambda *_{\mathrm{aux}} v^\prime)
  = \omega (v, v^\prime)\;.
\end{equation*}
Hence $\omega$ is ${\mathrm{aux}}$--invariant on the zero section of
$\nu_k$ and we do not change it there by averaging.  It follows that
there is a small neighbourhood of the zero section, where
$\omega^\prime$ will be symplectic.

\begin{propo}\label{properties of phi-action}
  There exists a Hamiltonian $\S^1$--action $\phi$ in a neighbourhood
  of $W_k$ with the following properties:
  \begin{itemize}
  \item[(i)] $\phi$ leaves $W_k$ pointwise fixed and the stabilizer of
    points not lying on $W_k$ is a proper subgroup of $\Z_k$;
  \item[(ii)] $\phi$ commutes with the $\beta$--action;
  \item[(iii)] the two actions coincide for elements in $\Z_k \le
    \S^1$, and if $\lambda*_\phi x = \sigma*_\beta x$, then
    $\sigma\in\Z_k\le \S^1$;
  \item[(iv)] the Hamiltonian function for the $\phi$--action is
    $\beta$--invariant and constant along $W_k$. Furthermore, it is a
    Morse--Bott function with critical set $W_k$ of index zero.
 \end{itemize}
\end{propo}
\begin{proof}
  We will obtain the $\phi$--action by a deformation of the existing
  auxiliary action.  By a standard application of the equivariant
  Moser trick, there exists a $\beta$--equivariant isotopy of a
  neighbourhood of $W_k$ that deforms $\omega$ into $\omega^\prime$.
  The only difficulty comes from the fact that we are dealing with an
  open set, so that the flow of the Moser vector field $X_t$ does not
  need to exist up to time~$1$.  But since $X_t$ vanishes on the zero
  section $W_k$, after shrinking $W_k$ and the radius of the tubular
  neighbourhood, we get a smaller set where the isotopy can be defined
  as the flow of $X_t$ for all $t\in [0,1]$. The inverse of this
  isotopy deforms the auxiliary action into the required action
  $\phi$.
  
  Since the flow of the Moser vector field leaves $W_k$ invariant and
  is $\beta$--equivariant, one can show that properties (i)--(iii)
  hold for the auxiliary action together with $\beta$ and this is
  equivalent to the corresponding statementes for $\phi$ and $\beta$.

  Since $\omega$ is $\phi$--invariant, one has that $di_{X_\phi}\omega
  = \lie{X_\phi}\omega = 0$.  For the time being, let $U$ be any
  tubular neighbourhood of $W_k$, where $\phi$ is defined.  The closed
  $1$--form $i_{X_\phi}\omega$ represents a class in $H^1(U)$ which
  vanishes if we pull it back to the zero section $W_k$: Given that
  $H^1(U) \cong H^1(W_k)$, it follows that $i_{X_\phi}\omega$ is exact
  on $U$, i.e., there exists a function $\mu_\phi$ such that
  $i_{X_\phi}\omega = d\mu_\phi$.  The function $\mu_\phi$ is uniquely
  defined up to an additive constant (which we may choose such that
  $\mu_\phi \equiv 0$ on $W_k$) and is $\beta$--invariant.

  Recall that a Hamiltonian $\S^1$--function is always Morse-Bott (see
  \cite[Section~5.5]{McDuffSalamonIntro}).  The critical set coincides
  with the set of fixed points of the action, hence in our case with
  $W_k$.

  The index of a Morse-Bott function is invariant under
  diffeomorphisms.  Therefore in order to compute the index of
  $\mu_\phi$ at $W_k$, it suffices to compute the index of the
  Hamiltonian function $\mu_{\mathrm{aux}}$ of the auxiliary action
  with respect to $\omega^\prime$, because $\mu_\phi$ is the pullback
  of $\mu_{\mathrm{aux}}$ under the Moser-flow.

  Let $\partial_r$ be the radial vector field on $\nu_k$ given by
  \begin{equation*}
    \partial_r(v) = \left.\frac{d}{dt}\right|_{t=1} t\cdot v
  \end{equation*}
  for $v\in \nu_k$.  We will prove that $\mu_{\mathrm{aux}}$ strictly
  increases in radial direction, which shows that it has index zero.
  More precisely, we will prove that $\lie{\partial_r}
  \mu_{\mathrm{aux}} \ge 0$ with equality only at the zero section.
  By definition of $\mu_{\mathrm{aux}}$, one has $i_{\partial_r}
  d\mu_{\mathrm{aux}} = \omega^\prime (X_{\mathrm{aux}}, \partial_r)$,
  so it will suffice to show that there exists a neighbourhood of
  $W_k$ where $\omega^\prime (X_{\mathrm{aux}}, \partial_r) \ge 0$.

  With $\pi$ denoting the bundle projection $\nu_k\to W_k$, the
  vertical bundle $V(\nu_k)\leq T\nu_k$ of $\nu_k$ can be identified
  with the pullback bundle
  \begin{equation*}
    \pi^* \nu_k = \bigl\{ (v,w) \in \nu_k \times \nu_k \bigm| \,
    \pi(v) = \pi(w)\bigr\} \;.
  \end{equation*}
  The identification of $\pi^*(\nu_k)$ and $V(\nu_k)$ goes as follows
  \begin{equation*}
  \chi:\, \pi^*(\nu_k) \to V(\nu_k), \quad (v,w) \mapsto
    \left.\frac{d}{dt}\right|_{t=0} (v+ tw) \;.
  \end{equation*}
  Let $v\in \nu_k$, and write it as $v = v_1 + \dotsm + v_l$ with
  respect to the splitting $\nu_k = E_1 \oplus \dotsm \oplus E_l$.
  Then the vectors $X_{\mathrm{aux}}$ and $\partial_r$ are given by
  \begin{equation*}
    X_{\mathrm{aux}}(v) = \chi (v, a_1 j v_1 + \dotsm + a_l j v_l) \quad 
    \text{ and }\quad \partial_r(v) = \chi (v,v)
  \end{equation*}
  as elements of $V(\nu_k) \cong \pi^*(\nu_k)$.  
  Now assume $v = (x,0)$ lies in the zero section of $\nu_k$.  Then
  \begin{equation*}
    \omega^\prime \bigl(\chi(v,\sum a_i j w_i), \chi(v,w_1+\dotsm+w_l)\bigr) =
     \sum_{j=1}^l a_j\, \omega^\prime \bigl(j \chi(v,w_j),
    \chi(v,w_j)\bigr) > 0
  \end{equation*}
  if $w\ne 0$, since the eigenspaces $E_j$'s are $\omega$--orthogonal.
  By continuity this also holds for all $v$ in a neighbourhood of the
  zero section and all $w\ne 0$.  Hence in particular $\omega^\prime
  (X_{\mathrm{aux}}, \partial_r) > 0$ on $U - W_k$.
\end{proof}

\subsection{Surgery along the minimal stratum}

In this section we describe how to replace a neighbourhood of (one
component of) the minimal stratum $W_k$ by a smooth manifold, and
extend both the symplectic form $\omega$ and the Hamiltonian
$\S^1$--action $\beta$ to the resulting manifold, in such a way that
the singular points of the extended action are of order strictly
smaller than $k$. In order to achieve this, we perform a symplectic
cut with respect to the Hamiltonian $\S^1$--action $\phi$ constructed
in the previous section. The main reference for symplectic cuts and
hence for everything that follows is Lerman's original paper
\cite{Lerman_SymplecticCuts}.

As in the previous section, consider the minimal singular stratum
$W_k$, denote by $U$ a tubular neighbourhood in $W$, where the
$\phi$--action is defined, and take now the product $U\times \C$.  It
admits a first circle action, which is just the extension of the
original $\S^1$--action $\beta$ by the trivial action on the
$\C$--factor, namely, for $x\in U$, $z\in\C$
\begin{equation*}
  \lambda *_\beta (x,z) :=(\lambda *_\beta x, z) \;,
\end{equation*}
and we can define a second circle action on $U\times \C$ by
setting
\begin{equation*}
  \lambda *_\phi (x,z) = (\lambda *_\phi x, \lambda^{-k}z) \;.
\end{equation*}
These two actions commute and therefore we can combine them and define
a new $\S^1$--action
\begin{equation*}
  \lambda *_\tau (x,z) :=\lambda *_\phi \left(\lambda^{-1} *_\beta
    (x,z)\right) \;.
\end{equation*}
This $\tau$--action is not effective, because the $\phi$-- and the
$\beta$--action coincide for elements in $\Z_k$.  Hence consider the
short exact sequence
\begin{equation*}
  0\to \Z_k\to\S^1\to \hat\S^1 \to 0 \;,
\end{equation*}
with the homomorphism of the circle given by $\lambda\mapsto
\lambda^k$, and let $\hat{\S}^1$ act on $U\times \C$ by
$\sigma*_{\hat\tau} (x,z)=\lambda*_\tau (x,z)$ for some $\lambda\in
\S^1$ such that $\lambda^k=\sigma$.  This new action, which we denote
by $\hat{\tau}$, is not only effective but by
Proposition~\ref{properties of phi-action}.(iii), even free and the
quotient $(U\times \C)/\hat{\tau}$ is a smooth manifold.  It still
carries an $\S^1$--action induced by the $\phi$--action on $U\times
\C$, and this is well defined because $\phi$ commutes with~$\hat\tau$.

\begin{figure}[htbp]
  \begin{center}
    \input{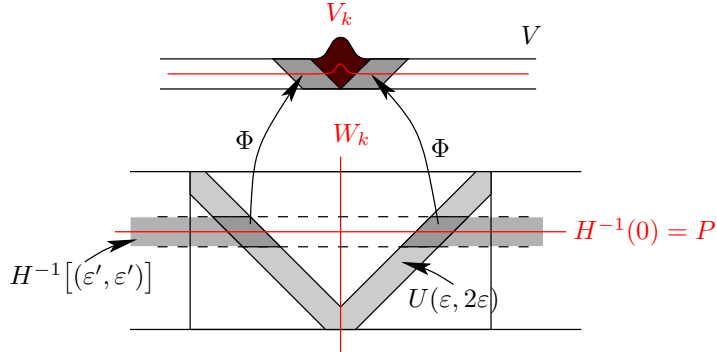}
    \caption{Remove the $U(\epsilon)$--neighbourhood of $W_k$, and
      glue in the $V(\epsilon)$--patch, identifying along
      $U(\epsilon,2\epsilon)$ via the diffeomorphism $\Phi$.}
    \label{fig: gluing patch}
  \end{center}
\end{figure}

We define a symplectic form $\Omega = (\omega, -i\,dz\wedge d\bar{z})$
on $U\times \C$, which is invariant with respect to the
$\hat\tau$--action.  By construction, the infinitesimal generator of
this action can be written as $X_{\hat\tau} = -X_{\beta}+X_{\phi}$.
The Hamiltonian for the ${\hat\tau}$--action is given by
\begin{equation*}
  H_{\hat\tau}(x,z) = \mu_\phi(x) - k\abs{z}^2 - H(x)\;,
\end{equation*}
and if we now do symplectic reduction at some level $\epsilon$ we get
the quotient $H_{\hat\tau}^{-1}(\epsilon)/\hat{\tau}$ which, with the
structure induced by $\Omega$ and $\phi$, is a smooth Hamiltonian
$\S^1$--manifold with Hamiltonian function $H_\phi[(x,z)] =
\mu_\phi(x) - k\abs{z}^2 = H(x) + \epsilon$.  Notice that $H_{\hat
  \tau}^{-1}(\epsilon)$ can be written as the disjoint union of two
$\hat{\tau}$--invariant manifolds
\begin{equation*}
  H_{\hat\tau}^{-1}(\epsilon)= \left\{(x,z) \left|\:
      \mu_\phi(x) - H(x) > \epsilon,\ 
      \abs{z}^2 = \frac{\mu_\phi(x)-H(x)-\epsilon}{k}\right.
  \right\}\sqcup
  \Bigl\{(x,0) \Bigm|\,\mu_\phi(x) - H(x) =\epsilon \Bigr\} \;.
\end{equation*}
Since the level sets of $\mu_\phi$ and $H$ are transverse, the
restriction of $\mu_\phi$ to $P = H^{-1}(0)$ is still a Morse-Bott
function of index zero.  Hence we can choose $\delta > 0$ such that
$\mu_\phi^{-1}(\delta) \cap P$ is contained in the interior of $U$ and
has the structure of a sphere bundle over $P_k = P\cap W_k$.  Choose
$\epsilon^\prime > 0$ such that
\begin{equation*}
  \bigl(\mu_\phi - H\bigr)^{-1} (\delta/2, \delta) \cap P_t
\end{equation*}
is non-empty, does not intersect $W_k$, and is contained in the
interior of $U$ for all $\abs{t} < \epsilon^\prime$ (see
Fig.~\ref{fig: gluing patch}).  Assume from now on that all subsets in
$W$ are restricted to $H^{-1}\bigl((-\epsilon^\prime,
\epsilon^\prime)\bigr)$.  Then
\begin{equation*}
  U(\delta/2,\delta) :=  \bigl(\mu_\phi - H\bigr)^{-1} (\delta/2, \delta)
\end{equation*}
is diffeomorphic to a spherical shell bundle over $W_k$.  Choose now
$\epsilon = \delta/2$, set $V := H_{\hat\tau}^{-1}(\epsilon) /
\hat\tau$, and consider the map
\begin{equation*}
  \Phi:\, U(\epsilon, 2\epsilon) \to V, \qquad x\mapsto 
  \left[x, \sqrt{\frac{\mu_\phi(x)- H(x) -\epsilon}{k}}\right] \;.
\end{equation*}
This is a diffeomorphism (onto its image), equivariant with respect to
the $\beta$--action on $U(\epsilon, 2\epsilon)$ and the $\phi$--action
on $V$.  Its inverse can be constructed as follows: given $[x,z]$ with
$z\neq 0$, we first represent the same class by an element $(x', z')$
such that $z'$ is a real positive number, and then define
\begin{equation*}
  \Phi^{-1}([x,z]):= x' \;.
\end{equation*}
Moreover, since $\Phi$ factors through a map $ U(\epsilon, 2\epsilon)
\to H_{\hat\tau}^{-1}(\epsilon)$ which is the identity in the first
component and a real function in the second one, we have
\begin{equation*}
  \Phi^*(\omega, -i\,dz\wedge d\bar{z})=\omega \;,
\end{equation*}
hence $\Phi$ gives in fact an equivariant symplectic identification of
$U(\epsilon, 2\epsilon)$ with its image under $\Phi$.  More precisely
we have
\begin{equation*}
  \Phi\bigl(U(\epsilon, 2\epsilon)\bigr)  =\bigl\{[x,z]\in V\,\bigm|\:
  \epsilon<\mu_\phi(x) - H(x) < 2\epsilon \bigr\}  \;.
\end{equation*}
We can now remove a tubular $\epsilon$--neighbourhood
$\overline{U(\epsilon)}$ of $W_k$ in $W$ and glue in the smooth
manifold
\begin{equation*}
  V(\epsilon):=
  \left.\left\{(x,z)\,\left|\:\epsilon\leq\mu_\phi(x) - H(x) < 2\epsilon,
        \,
        \abs{z}^2=\frac{\mu_\phi(x)-H(x)-\epsilon}{k} \right.\right\}
  \right/\hat{\tau}
\end{equation*}
along the open ``collar'' $U(\epsilon, 2\epsilon)$, using the map
$\Phi$.  In this way we define the new manifold
\begin{equation*}
  \widetilde W =\bigl(W-\overline{U(\epsilon)}\bigr)\cup_{\Phi} 
  V(\epsilon) \;.
\end{equation*}
Since $\Phi$ is equivariant, the $\beta$--action on $W -
\overline{U(\epsilon)}$ and the $\phi$--action on $V(\epsilon)$ fit
together to give a circle action $\widetilde\beta$ on $\widetilde W$,
which by construction coincides with $\beta$ outside a
$2\epsilon$--neighbourhood of $W_k$.  Moreover, $\Phi$ identifies the
given symplectic forms on the two sides of the gluing, so $\widetilde
W$ also admits a symplectic form $\widetilde{\omega}$ with the
property that $\widetilde{\omega}=\omega$ on $W- U(2\epsilon)$.  With
the action $\widetilde\beta$ and the symplectic form
$\widetilde{\omega}$ just defined, $\widetilde W$ is a Hamiltonian
$\S^1$--manifold.  The Hamiltonian function $\widetilde H$ for
$\widetilde\beta$ is given by $H$ on $W - U(2\epsilon)$ and by $H_\phi
- \epsilon$ on $V(\epsilon)$.

We need to analyse the singular points of the
$\widetilde{\beta}$--action on the ``patch'' $V(\epsilon)$.  They
satisfy the relation $\lambda*_{\phi}[x,z]=[x,z]$ for some
$\lambda\in\S^1$ and this in turn means that there exist $\kappa\in
\hat{\S}^1$ and $\sigma\in \S^1$, $\sigma^k=\kappa$ such that
\begin{equation*}
  \lambda*_{\phi}(x,z)=\kappa*_{\hat{\tau}}(x,z)=\sigma*_{\phi}
  \sigma^{-1}*_{\beta}(x,z) \;.
\end{equation*}
In particular, $\lambda*_{\phi}x =\sigma*_{\phi} \sigma^{-1}*_{\beta}
x$, and by Proposition~\ref{properties of phi-action}.(iii) this
identity can only hold if $\sigma\in \Z_k$.  Hence $\kappa = 1$ and
singular points are characterised by $\lambda*_{\phi}(x,z)=(x,z)$.
Since $[x,z]\in V(\epsilon)$ implies that $x$ does not lie in $W_k$,
the isotropy groups are proper subgroups of $\Z_k$, see
Proposition~\ref{properties of phi-action}.(i).  

We have thus proved:

\begin{propo}
Let $W_k$ be the minimal singular stratum of a Hamiltonian 
$\S^1$--manifold $W$. There exists a new symplectic manifold 
$\widetilde{W}$ which is symplectomorphic to $W$ outside an 
arbitrarily small neighbourhood of $W_k$ and admits a Hamiltonian 
$\S^1$--action that coincides with the given action away from $W_k$ 
and only has singularities of order strictly smaller than $k$.
\end{propo}

In fact we have also shown that $0$ is still a regular value of the
moment map of the circle action on $\widetilde{W}$ and if we consider
the symplectic reduced space at this level, namely
\begin{equation*}
  \widetilde M := \widetilde H^{-1}(0)/\widetilde\beta \;,
\end{equation*}
we see that this symplectic orbifold has singularities of strictly
lower order than those of $M = H^{-1}(0)/\beta$.

Moreover, there exists a map $f:\,\widetilde{M}\to M$, which is a
symplectic orbifold isomorphism outside an arbitrarily small
neighbourhood of $M_k = P_k/\beta$ (and in fact coincides with the
identity map outside a slightly larger neighbourhood).  We shall
describe how to define $f$.  On $(P-U(2\epsilon))/\beta$ it is simply
the identity.  In order to define it on $\bigl(V(\epsilon)\cap
H^{-1}_\phi(\epsilon)\bigr)/\phi$ a little more work is required.
First of all, denote by $V_k$ the quotient
\begin{equation*}
  \bigl\{(x,0)\in
  H_{\hat\tau}^{-1}(\epsilon)\,\bigm|\,\mu_\phi(x) - H(x) =
  \epsilon \bigr\}/\hat\tau \;.
\end{equation*}
Then the inverse of the gluing map $\Phi$ restricts to a
diffeomorphism $\Phi^{-1}:\, (V-V_k)\cap H^{-1}_\phi(\epsilon) \to
P\cap U(\epsilon,2\epsilon)$.  Since $\Phi$ is equivariant with
respect to the $\phi$-- and $\beta$--actions, this descends to a
symplectic orbifold isomorphism on the quotients.  Let $h:\,
U(\epsilon,2\epsilon) \cap P \to (U(2\epsilon) - W_k) \cap P$ be a
$\beta$-- and $\phi$--equivariant diffeomorphism (recall that $\phi$
is the action on $W$ constructed in Section~\ref{sec: construction of
  auxiliary circle action}), that is the identity in a neighbourhood
of the outer boundary of $U(2\epsilon)$ and extends to a smooth map
from $U(2\epsilon) \cap P$ to itself which maps $U(\epsilon)\cap P$ to
$P_k$.

Then we can define
\begin{align*}
  f: \quad & \bigl(V\cap H^{-1}_\phi(\epsilon)\bigr) /\phi \to
  \bigl(P\cap U(2\epsilon)\bigr) / \beta \\
  & [x,z] \mod \phi \mapsto
  \begin{cases}
    h\circ \Phi^{-1} ([x,z]) \mod \beta & \text{ if $z\ne 0$,} \\
    h(x) \mod \beta & \text{ if $z=0$.}
  \end{cases}
\end{align*}
Because of the boundary conditions on $h$, the map $f$ extends on the
outer side to the identity map.  To see that $f$ is continuous in a
neighbourhood of $\bigl(V_k\cap H^{-1}_\phi(\epsilon)\bigr) / \phi$,
one has to show that for any sequence $[x_k,z_k]/\phi \subset
\bigl(V\cap H^{-1}_\phi(\epsilon)\bigr)/\phi$ that converges to some
element $[x_0,0]/\phi$, it follows that $f\bigl([x_k,z_k]/\phi\bigr)$
converges to $f\bigl([x_0,0]/\phi\bigr)$.  We can find representatives
$(x_k^\prime,z_k^\prime)$ for the sequence that converge to
$(x_0^\prime,0)$ and such that $z_k^\prime$ is a non-negative real
number, and hence $f\bigl([x_k^\prime,z_k^\prime]/\phi\bigr) =
h(x_k^\prime) /\beta$ converges to $h(x_0^\prime)/\beta =
f([x_0^\prime,0]/\phi) = f([x_0,0]/\phi)$.

\subsection{Proof of Theorem \ref{main}}
If we inductively repeat the above surgery procedure along the minimal
stratum, we will eventually obtain a symplectic manifold
$\widetilde{W}$ with a free Hamiltonian $\S^1$--action
$\widetilde{\beta}$.  The symplectic reduced space of this action at
the level $0$ will be a smooth symplectic manifold $\widetilde{M}$ and
it will be equipped with a function $f:\widetilde{M}\to M$, which
factors through all the previous steps of the resolution (that is,
reduced spaces with orbifold singularities of decreasing order), and
is a symplectic diffeomorphism outside an arbitrarily small
neighbourhood of the set of orbifold singularities of~$M$.


\section{Applications}

\subsection{Isolated cyclic orbifold singularities}

Any isolated cyclic orbifold singularity can be represented as
$\C^n/\Z_k$ with symplectic form $\omega = i\sum_{j=1}^n dz_j\wedge
d\bar z_j$, where the generator $\xi = e^{2\pi i/k}$ of $\Z_k$ acts by
\begin{equation*}
  \xi\cdot \z = (e^{2\pi i/k} z_1,
  e^{2\pi i a_2/k} z_2,\dotsc,
  e^{2\pi i a_n/k} z_n) \;,
\end{equation*}
and $a_2,\dotsc,a_n\in\N$ are all coprime with~$k$.

To find a resolution of the singularity using Theorem~\ref{main}, just
note that $\C^n/\Z_k$ can be obtained by doing symplectic reduction on
the manifold $\C^*\times \C^n$ with $\S^1$--action given by
\begin{equation*}
  e^{i\phi}\cdot (\rho e^{i\theta},z_1,\dotsc,z_n) := (\rho e^{i(k\phi+\theta)},
  e^{i\phi} z_1, e^{i a_2\phi}z_2,\dotsc, e^{i a_n\phi}z_n) \;,
\end{equation*}
symplectic form $\omega + \rho\, d\rho\wedge d\theta +
\frac{1}{k}\,d\theta\wedge \sum_{j=1}^n a_j\,d\abs{z_j}^2$ and moment
map $H(\rho e^{i\theta},z_1,\dotsc,z_n)=\rho^2$.  Theorem~\ref{main}
then immediately implies the following generalisation of the
desingularisation result obtained in dimension~$4$:

\begin{coro}
  Every symplectic orbifold with only isolated cyclic singularities
  admits a symplectic resolution.
\end{coro}

Cavalcanti, Fernández, and Muñoz 
\cite{Cavalcanti_Fernandez_Munoz_NonFormal} recently introduced a 
method that allows one to find resolutions for all isolated symplectic 
orbifold singularities.
		
\subsection{Generalised Boothby-Wang fibrations are fillable}

\begin{defi}
  A \textbf{Boothby-Wang fibration} is a closed contact manifold
  $(P,\alpha)$ with a free $\S^1$--action which is given by the flow
  of the Reeb field $\Reeb$.  A \textbf{generalised Boothby-Wang
    fibration} is a closed contact manifold $(P,\alpha)$, where the
  Reeb field induces a semi-free $\S^1$--action.
\end{defi}

\begin{remark}\label{remark: pre-quantization and Boothby-Wang} A
  Boothby-Wang fibration $(P,\alpha)$ defines an $\S^1$--principal
  bundle over the manifold $B = P/\S^1$ with connection $1$--form
  $\alpha$.  The curvature form $\omega$ is the unique $2$--form on
  $B$ that satisfies $\pi^*\omega = d\alpha$.  The base manifold
  $(B,\omega)$ is a symplectic manifold and $\omega$ represents an
  integral cohomology class.

  Conversely, for any symplectic manifold $(B,\omega)$ with integral
  symplectic form, one can construct a Boothby-Wang fibration
  $(P,\alpha)$ over it, the so-called \textbf{pre-quantisation}.  This
  is the inverse of the previous construction.
\end{remark}

\begin{remark}
  A generalised Boothby-Wang fibration can be considered as the
  pre-quantisation of the symplectic \emph{orbifold} $(P/\S^1,
  \omega)$, and all the statements made in Remark~\ref{remark:
    pre-quantization and Boothby-Wang} can be translated to this
  setting.
\end{remark}

\begin{propo}
  A generalised Boothby-Wang fibration $(P,\alpha)$ has a natural
  convex filling by a symplectic orbifold.
\end{propo}
\begin{proof}
  These computations were obtained with the help of H.~Geiges.
  Consider the (complex) ``line bundle'' $L$ associated to $P$, i.e.\
  the bundle obtained from $P\times\C$ by identifying $(p,z)$ with
  $(e^{-i\phi}*p,e^{i\phi}\,z)$ for every $e^{i\phi}\in\S^1$.  The
  manifold $P$ embeds naturally via
  \begin{equation*}
    P\hookrightarrow L,\quad p \mapsto [p,1]\;.
  \end{equation*}
  The two forms
  \begin{equation*}
    \frac{1}{2}\,\left(\abs{z}^2\,\alpha + x\,dy-y\,dx\right) \quad
    \text{ and }\quad \frac{1}{2}\,d\alpha
  \end{equation*}
  on $P\times\C$ induce well-defined forms on $L$.  By adding the
  differential of the first form to the second one, we obtain a
  pre-symplectic form
  \begin{equation*}
    \omega := \frac{1}{2}\,d(\abs{z}^2)\wedge\alpha + dx\wedge dy +
    \frac{1+\abs{z}^2}{2}\,d\alpha \;,
  \end{equation*}
  because $2^n\,\omega^n = n\,(1+\abs{z}^2)^{n-1}\,(d\alpha)^{n-1}
  \wedge \left(d\abs{z}^2\wedge \alpha + 2dx\wedge dy\right)$ has only
  a one-dimensional kernel on $P\times\C$ that is generated by $-Z_P +
  x\,\partial_y - y\,\partial_x$.  It follows that
  $(P\times\C,\omega)$ is a pre-symplectic $\S^1$--manifold, and hence
  $L$ is a symplectic orbifold where all orbifold singularities sit
  along the zero section.
  
  Finally, the following field
  \begin{equation*}
    X := \frac{1+r^2}{2r}\,\partial_r = \frac{1+ x^2 + y^2}{2\,(x^2
      + y^2)}\,\left(x\,\partial_x + y\,\partial_y\right)
  \end{equation*}
  is a Liouville vector field for the manifold $(P,\alpha)$, and
  $(L,\omega)$ is hence a convex filling of~$P$.
\end{proof}

All the orbifold singularities lie in the interior of $L$.  By passing
to a symplectic resolution of $L$ whose existence is guaranteed by
Theorem~\ref{main}, we can obtain a smooth symplectic filling.
 
\begin{coro}\label{generalized Boothby-Wang are fillable} Generalised
  Boothby-Wang fibrations are symplectically fillable by a
  \emph{smooth} manifold.
\end{coro}

\begin{remark}
  \begin{enumerate}
  \item Popescu-Pampu recently proved the following conjecture of
    Biran: there exist prequantizations (for example on higher
    dimensional tori) that are not holomorphically fillable
    \cite{PompescuPampuFillability}. These manifolds, on the other
    hand, do have a strong symplectic filling, showing that the
    different types of fillability do not coincide (a result
    well-known in dimension~$3$).
  \item Massot showed that any contact structure on a
    $3$--dimensional Seifert manifold that is transverse to the fibres
    (but not necessarily invariant) has a weak symplectic filling
    \cite{MassotGeodesibleContact}.  To achieve this result he
    constructs first a filling by an orbifold, whose singularities can
    then be resolved by using for example the method presented in this
    paper.
  \end{enumerate}
\end{remark}

\bibliographystyle{amsalpha}

\bibliography{main}

\end{document}